\documentclass[11pt,reqno]{amsart}
\usepackage{amsthm, amssymb, latexsym, multicol, verbatim, enumerate, accents, float, hyperref, amsmath, amscd, youngtab}
\usepackage[usenames]{color}
\usepackage[utf8]{inputenc}

\advance\textwidth by 1.2in  \advance\oddsidemargin by -.6in \advance\evensidemargin by -.6in

\theoremstyle{plain}
\newtheorem*{prop}{Proposition}
\newtheorem{thm}{Theorem}
\newtheorem*{lem}{Lemma}
\newtheorem*{cor}{Corollary}

\theoremstyle{definition}
\newtheorem*{example}{Example}
\newtheorem*{defn}{Definition}
\newtheorem*{theom}{Theorem}
\newtheorem*{rem}{Remark}

\theoremstyle{remark}

\newcommand{\lie}[1]{\mathfrak{#1}}   \newcommand\bc{\mathbb C} \newcommand\bn{\mathbb N} \newcommand\bz{\mathbb Z}   
   \def\gr{\operatorname{gr}}

\newcounter{cnt}
 \makeatletter
\def\mydggeometry{\makeatletter\dg@YGRID=1\dg@XGRID=20\unitlength=0.003pt\makeatother}
\makeatother \theoremstyle{remark}

\numberwithin{equation}{section}
\makeatletter
\def\section{\def\@secnumfont{\mdseries}\@startsection{section}{1}%
  \z@{.7\linespacing\@plus\linespacing}{.5\linespacing}%
  {\normalfont\scshape\centering}}
\def\subsection{\def\@secnumfont{\bfseries}\@startsection{subsection}{2}%
  {\parindent}{.5\linespacing\@plus.7\linespacing}{-.5em}%
  {\normalfont\bfseries}}
\makeatother

\begin{document}
\title[PBW degenerations of Lie superalgebras]{PBW degenerations of Lie superalgebras and their typical representations}

\author{Ghislain Fourier}
\address{RWTH Aachen University, Pontdriesch 10-16, 52062 Aachen}
\email{fourier@mathb.rwth-aachen.de}
\thanks{}
\author{Deniz Kus}
\address{Ruhr-University Bochum, Universitätsstr. 150, 44780 Bochum}
\email{deniz.kus@rub.de}
\thanks{}

\subjclass[2010]{}
\begin{abstract}
We introduce the PBW degeneration for basic classical Lie superalgebras and construct for all type I, $\mathfrak{osp}(1,2n)$ and exceptional Lie superalgebras new monomial bases. These bases are parametrized by lattice points in convex lattice polytopes, sharing useful properties such as the integer decomposition property. This paper is the first step towards extending the framework of PBW degenerations to the Lie superalgebra setting.
\end{abstract}
\maketitle \thispagestyle{empty}
\section{Introduction}
The framework on PBW degenerations in Lie theory started around ten years ago with several works of Evgeny Feigin (\cite{F2009, FFL11, FFL12}). Roughly speaking, the main idea is to \textit{degenerate} a Lie algebra into an abelian Lie algebra. For example, on the level of universal enveloping algebras the PBW degree of monomials can be used to define a filtration such that the corresponding associated graded algebra is isomorphic to an ordinary polynomial ring. 
There is an induced filtration on any cyclic module and the induced associated graded space is currently one main object of research. This has been developed for finite-dimensional, simple complex Lie algebras and their finite-dimensional simple modules. Here, the associated graded space is a quotient of the polynomial ring, the commutative analogue of the well--known description of a finite-dimensional simple module as a quotient of the universal enveloping algebra.\\ This description raises a couple of interesting questions. Can we determine generators for the defining ideal? Can we provide a monomial basis for the associated graded space? In the classical, non-degenerate setup, these answers are known for a long time. In the degenerate setup, there are positive answers in type $A$ \cite{FFL11}, type $C$ \cite{FFL12} and type $G$ \cite{G15}. Beyond these cases little is known, there hasn't been any new result on infinite series for the last seven years. \par
In the present paper, we will extend the framework of PBW degenerations to Lie superalgebras. The notion of PBW degree of monomials in the universal enveloping algebra is natural in this context, the elements of a fixed basis will have degree one. The associated graded algebra is then a product of a symmetric algebra with an exterior algebra. Again, one has an induced filtration on cyclic modules and the natural questions of finding generators for the defining ideal and monomial bases show up in this context too.

We restrict ourselves to finite-dimensional typical representations $V(\lambda)$ of a basic classical Lie superalgebra $\lie g$; for the necessary choices of positive roots and Borel subalgebras we refer the reader to Section~\ref{section2}. For the moment we assume an appropriate triangular decomposition $\lie g = \lie n^+ \oplus \lie h \oplus \lie n^-$  and set $\lie n^{\pm}_{\bar{i}}=\lie g_{\bar{i}}\cap \lie n^{\pm}$. We consider the PBW degenerate Lie superalgebra $\lie n^{-,a}$ and the induced PBW degenerate module $V^a(\lambda)$. 
Recall that $\lie n^{-,a}$ is isomorphic to $\lie n^{-}$ as a vector space with trivial Lie superbracket and hence the universal enveloping algebra $\mathbf{U}(\lie n^{-,a})$ can be identified with $S(\lie n^{-}_{\bar{0}})\otimes \Lambda(\lie n_{\bar{1}}^-)$. The PBW filtration is compatible with the triangular decomposition in the sense that $\mathbf{U}(\lie n^{-,a})$ and $V^a(\lambda)$ are natural $\mathbf{U}(\lie n^+)$--modules. Let us first consider the type I case, for technical reason we do not consider $A(n,n)$ here but our results can be extended without any difficulties for the central extension of $A(n,n)$. 
In this setup, the simple ideals in the underlying reductive Lie algebra $\lie g_{\bar{0}}$ are of type $A$ or $C$ and hence a monomial basis for the PBW degenerate module of $\lie g_{\bar{0}}$ are described using the so-called FFLV-polytopes, introduced in \cite{FFL11, FFL12}. For fixed $\lambda$, we denote this polytope by $P_{\lie g_{\bar{0}}}(\lambda)$ and its lattice points by $S_{\lie g_{\bar{0}}}(\lambda)$. 
We extend the results from classical Lie algebras to provide a complete answer for Lie superalgebras of type I. 
\begin{theom}
Let $\lambda, \mu \in P^+$ be dominant integral typical weights and $d$ the number of positive odd roots. We set $P(\lambda) := P_{\lie g_{\bar{0}}}(\lambda) \times [0,1]^{d}$ and $S(\lambda)$ denotes the lattice points in $P(\lambda)$. 
\begin{enumerate}
\item We have $P(\lambda + \mu) \subseteq P(\lambda) + P(\mu)$ and $S(\lambda + \mu) \subseteq S(\lambda) + S(\mu)$.
\item The lattice points $S(\lambda) \subseteq P(\lambda)$ parametrize a monomial basis for $V^a(\lambda)$.
\item We have $V^a(\lambda + \mu) \subseteq V^a(\lambda) \otimes  V^a(\mu)$.
\item We have an isomorphism $V^a(\lambda) \cong \mathbf{U}(\lie n^{-,a})/\mathbf{I}(\lambda)$, where $\mathbf{I}(\lambda)$ is generated by the set 
$$\Big\{ \mathbf{U}(\lie n_{\bar{0}}^+) \circ  x_{-\alpha}^{2\frac{(\lambda,\alpha)}{(\alpha,\alpha)} + 1} : \alpha \in R_{\bar{0}}^+ \Big\} \subseteq \mathbf{U}(\lie n^{-,a}).$$
\end{enumerate}
\end{theom}

For Lie superalgebras of type II, the answer is more complicated, as $V(\lambda)$ is a proper quotient of the induced module (see Remark~\ref{grad}). Our first result reduces the problem of finding a monomial basis for the associated graded space to the computation of monomial bases for $\gr V(\lambda)$, where the highest weight $\lambda$ is supported only on a unique simple ideal of $\lie g_{\bar{0}}$, see Proposition~\ref{reduction}. This reduction is used for example in Section~\ref{section-appendix} to construct a bases in the exceptional cases.
We explain our results in the case of $\lie{osp}(1,2n)$, which will be the first new infinite series to be solved since \cite{FFL12}. Inspired by the definition of symplectic Dyck path, we introduce the notion of orthosymplectic Dyck paths (see Definition~\ref{dyckpath}) and define for $\lambda \in P^+$ a convex lattice polytope $P_{\lie{osp}}(\lambda)$ and denote by $S_{\lie{osp}}(\lambda)$ the corresponding lattice points. In contrast to the type $A$ and $C$ simple Lie algebras, this is not a marked chain polytope (see \cite{ABS11}). Nevertheless it shares certain useful properties with those. Our main result is the following.
\begin{theom}
All the statements of the previous theorem are true for $\lie{osp}(1,2n)$ using $P_{\lie{osp}}(\lambda)$ and $S_{\lie{osp}}(\lambda)$. For the precise description of the ideal $\mathbf{I}(\lambda)$, see Theorem~\ref{thm1}.
\end{theom}
Recall, that in the classical Lie algebra cases $A$ and $C$, one has $S(\lambda + \mu) = S(\lambda) + S(\mu)$. In the super setting, this can't be an equality, as the power of positive odd root vectors is limited to $1$. Define
\[
\Sigma = \bigcup_{\lambda \in P^+ } S(\lambda) \times \{ \lambda \}, \; \;  \Sigma_{\lie{osp}} = \bigcup_{\lambda \in P^+ } S_{\lie{osp}}(\lambda) \times \{ \lambda \}, 
\] 
then $P(\lambda)$ (resp. $P_{\lie{osp}}(\lambda)$) is a slice of the polyhderal cone defined by the convex hull of $\Sigma$ (resp. $\Sigma_{\lie{osp}}$). We have shown in particular, that the semi group $\Sigma$ (resp. $\Sigma_{\lie{osp}}$) is finitely generated by the lattice points for fundamental weights $S(\varpi_i) \times \{ \varpi_i\}$ (resp. $S_{\lie{osp}}(\varpi_i) \times \{ \varpi_i\}$).\par
In the remaining infinite Lie superalgebra series of type II, we have statements reducing the problem of finding monomial bases for $V^{a}(\lambda)$ to the computation of monomial bases for PBW degenerate modules for simple Lie algebras of type $B$. Since these monomial bases are not yet described in full generality (see \cite{BK19} for partial results), we omit our reduction statements for now and refer to future publications. \par

Nevertheless, we are able to provide similar results for the exceptional types $F(4), G(3)$ and $D(2,1;\alpha)$. Here the results are less uniform as the semi groups of typical highest weights are not finitely generated. We provide for all typical representations a convex lattice polytope whose lattice points parametrize a monomial bases of the PBW degenerate typical modules (see Section~\ref{section-appendix} for more details). The methods in this section are different from the ones in Section~\ref{section3} and Section~\ref{section5}. We use the reduction procedure and count lattice points in convex simplex-like polytopes.\par

From the first papers on PBW degenerations for classical Lie algebras, there is a large variety of various applications, influences and appearances. We discuss briefly where to go from here. Feigin introduced the PBW degenerate flag variety, considering the orbit of the degenerate $SL_n$ on the PBW degenerate module (see \cite{Fei11, Fei12}). In here the links to combinatorics and quiver Grassmannian, linear degenerations of flag varieties (see \cite{CIFFFR17,CIFR13,CIFR12}) have been already provided.  
The combinatorics of the monomial bases is providing the link to the theory of crystal bases (see \cite{K12,K13}), toric degenerations of flag varieties (see \cite{FaFL17}), Newton-Okounkov bodies (see \cite{FFL17}) and discrete geometry (see \cite{ABS11}). We would like to think of the present paper as being the starting point of a similar study for super flag varieties (see for example \cite{P90}).
Does a PBW degenerate super flag variety would also have such an impact? How to interpret the combinatorics of the monomial bases? Would a notion of super marked chain polytopes be reasonable here? We restrict ourselves for this paper to the PBW filtration and monomial bases and postpone these questions to a forthcoming publication.\par

The paper is structured as follows: in Section~\ref{section2} we recall the most important definitions for Lie superalgebras and their typical representations, as well as the PBW filtration.  In Section~\ref{section3}, we consider the type I case, focusing on $A(m,n)$ and $C(n)$. In Section~\ref{section5}  we analyze the $\mathfrak{osp}(1,2n)$-case, providing the monomial bases and essential polytope. In the Appendix~\ref{section-appendix}, we consider the exceptional cases $F(4), G(3)$ and $D(2,1;  \alpha)$.\\

\textbf{Acknowledgements}: \textit{Part of this work was done when both authors visited the Catholic University of America for the conference Interactions of quantum affine algebras with cluster algebras, current algebras and categorification. They thank UCA for the superb working conditions and Prasad Senesi for the organization of the conference. G.F. would like to dedicate this paper to the SFB/TR 195 on ''Symbolic tools in mathematics and their application''.}
%

\section{Lie superalgebras, their representations and the PBW filtration}\label{section2}
\subsection{}We denote the set of complex numbers by $\mathbb{C}$ and, respectively, the set of integers, non--negative integers, and positive integers  by $\mathbb{Z}$, $\mathbb{Z}_+$, and $\mathbb{N}$. Unless otherwise stated, all the vector spaces considered in this paper are $\mathbb{C}$--vector spaces and $\otimes$ stands for $\otimes_\mathbb{C}$.
\subsection{}Finite-dimensional simple complex Lie superalgebras were classified by Kac \cite{K77} and are divided into two groups: the classical Lie superalgebras (the even part is a reductive Lie algebra) and the Cartan type Lie superalgebras. In this paper we consider basic classical Lie superalgebras, i.e. classical Lie superalgebras on which there exists a non-degenerate, invariant, even bilinear form. A basic classical Lie superalgebra $\lie g=\lie g_{\bar{0}}\oplus \lie g_{\bar{1}}$ is said to be of type II if $\lie g_{\bar{1}}$ is an irreducible $\lie g_{\bar{0}}$--module via the adjoint action, and is said to be of type I if $\lie g_{\bar{1}}$ is a direct sum of two irreducible $\lie g_{\bar{0}}$-modules. The following table gives all possible types of basic classical Lie superalgebras, which are not simple Lie algebras:
\begin{table}[H]
\begin{tabular}{lclcl}

$\lie g$ & $\lie g_{\bar{0}}$ \\
\hline
    $\Big. A(m,n)$, \ $m>n \geq 0$ & $A_m \oplus A_n \oplus \mathbb{C}$ &type I \\
    $\Big. A(n,n)$, \ $n \geq 1$ & $A_n \oplus A_n$ &type I \\
    $\Big. B(m,n)$, \ $m \geq 0$, $n \geq 1$ & $B_m \oplus C_n$ &  type II \\
    $\Big. C(n)$, \ $n \geq 2$ & $C_{n-1} \oplus \mathbb{C}$ & type I \\
    $\Big.D(m,n)$, \ $m \geq 2$, $n \geq 1$ & $D_m \oplus C_n$ &type II \\
    $\Big. D(2,1;\alpha)$, \ $\alpha \neq 0,-1$ & $A_1 \oplus A_1 \oplus A_1$ &type II\\
    $\Big.F(4)$ & $A_1 \oplus B_3$ & type II \\
    $\Big.G(3)$ & $A_1 \oplus G_2$ & type II \\

\end{tabular}
\bigskip
\end{table}

In terms of matrices we have (all entries are in the field of complex numbers) $A(m,n)=\mathfrak{sl}(m+1,n+1)$, $A(n,n)=\mathfrak{psl}(n+1,n+1)$, $B(m,n)=\mathfrak{osp}(2m+1,n+1)$, $C(n)=\mathfrak{osp}_(2,2n-1)$, $D(m,n)=\mathfrak{osp}(2m,2n)$. \textit{For technical reasons we want to leave out the Lie superalgebra $A(n,n)$ in the rest of this paper, but we can extend without any difficulties our results for the central extension $\mathfrak{sl}(n+1,n+1)$ of $A(n,n)$.}
\subsection{}We fix a Cartan subalgebra $\lie h\subseteq \lie g$, which is by definition a Cartan subalgebra of $\lie g_{\bar{0}}$. For $\alpha\in \lie h^{*}$ let
$$\lie g_{\alpha}=\big\{x\in \lie g : [h,x]=\alpha(h)x,\ \forall  h\in \lie h\big\}$$
the root space associated to $\alpha$ and $R=\big\{\alpha\in \lie h^{*}\backslash\{0\}: \lie g_{\alpha}\neq 0\big\}$ be the root system of $\lie g.$ We obtain
\begin{equation}\label{roots}\lie g=\lie h\oplus \bigoplus_{\alpha\in R}\lie g_{\alpha},\end{equation}
where each root space $\lie g_{\alpha}$ in \eqref{roots} is one--dimensional \cite[Proposition 1.3]{K78}. 
We define the even and odd roots to be
$$R_{\bar{0}}=\big\{\alpha\in R \mid \lie g_{\alpha}\cap \lie g_{\bar{0}}\neq 0\big\},\ R_{\bar{1}}=\big\{\alpha\in R \mid \lie g_{\alpha}\cap \lie g_{\bar{1}}\neq 0\big\}.$$
Let $E$ be the real vector space spanned by $R$ equipped with a total ordering $\succ$ compatible with the real vector space structure. We denote by $R^{+}=\{\alpha\in R \mid \alpha \succ 0\}$ and $R^{-}=\{\alpha\in R \mid \alpha \prec 0\}$ respectively the set of positive roots and negative roots respectively. We fix a subset $\Delta=\{\alpha_1,\dots,\alpha_r\}\subseteq R^+$ of simple roots, which by definition means that $\alpha\in \Delta$ cannot be written as a sum of two positive roots. We denote by $I=\{1,\dots,r\}$ the corresponding index set. Let $\rho_{\bar{0}}$ (respectively $\rho_{\bar{1}}$) be the half--sum of all the even (respectively odd) positive roots and set $\rho=\rho_{\bar{0}}-\rho_{\bar{1}}$.
We have a triangular decomposition 
$$\lie g=\lie n^-\oplus \lie h \oplus \lie n^+,\mbox{ where } \lie n^{\pm}=\bigoplus_{\alpha\in R^{\pm}} \lie g_{\alpha}.$$
The subalgebra $\lie b=\lie h\oplus \lie n^+$ is called the Borel subalgebra of $\lie g$ corresponding to the positive system $R^+$. We emphasize that unlike in the setting of semi--simple Lie algebras, Borel subalgebras need not to be conjugate. Most of the theory heavily depends on the choice of a simple system, but we will see later that there is a canonical choice which is called distinguished simple system. We set 
$$\lie n^{\pm}_{\bar{i}}=\lie n^{\pm}\cap \lie g_{\bar{i}},\ \ 0\leq i \leq 1.$$
\begin{example}\label{ex1}We discuss the properties of the Lie superalgebra 
$$\lie g=\mathfrak{sl}(m,n)=\left\{\begin{pmatrix}
A & B \\
C & D  \\
\end{pmatrix}: A\in \mathbb{C}^{m\times m}, B\in \mathbb{C}^{m\times n}, C\in \mathbb{C}^{n\times m}, D\in \mathbb{C}^{n\times n}, str(A)=0\right\},$$
where $str(A):=tr(A)-tr(D)$ denotes the supertrace of $A$. 
The subspace consisting of matrices where $B=C=0$ determines a reductive Lie algebra isomorphic to $\lie g_{\bar{0}}$ and since $\lie g_{\bar{1}}$ decomposes into two irreducible $\lie g_{\bar{0}}$-modules we have that $\lie g$ is a basic classical Lie superalgebra of type I. The non-degenerate bilinear form is given by the formula
$$(\cdot,\cdot): \lie g\times \lie g\rightarrow \mathbb{C}, (A,B)\mapsto str(AB).$$
The diagonal matrices in $\lie g$ form a Cartan algebra and the corresponding roots are given by
$$R_{\bar{0}}=\{\epsilon_i-\epsilon_j,\ \delta_r-\delta_s : 1\leq i\neq j\leq m,\ 1\leq r\neq s\leq n\},\ \ R_{\bar{1}}=\{\pm(\epsilon_i-\delta_j): 1\leq i\leq m,\ 1\leq j\leq n\},$$
where for $h=\text{diag}(a_1,\dots,a_m,b_1,\dots,b_n)$, $\epsilon_i(h)=a_i$ and $\delta_i(h)=b_{i}$. There are several choices for simple systems, but the most canonical one is the following with exactly one odd root:
\begin{equation}\label{ss1}\Delta=\{\delta_1-\delta_2, \dots, \delta_{n-1}-\delta_n, \delta_n-\epsilon_1, 
\epsilon_1-\epsilon_2, \dots, \epsilon_{m-1}-\epsilon_{m}\}.\end{equation}
\end{example}
\subsection{}\label{section24}
A postive root system is called distinguished if the corresponding system of simple roots contains exactly one odd root; for example the set of simple roots in \eqref{ss1} is distinguished. \textit{From now on we fix a distinguished positive root system for $\lie g$ with Cartan matrix $A=(a_{i,j})_{i,j\in I}$ whose Dynkin diagram $S$ is given as in \cite[Table 1]{K78}. We denote by $s$ the unique node such that $\alpha_s$ is odd.} The following table gives an explicit description of the distinguished simple system $\Delta$.

\begin{table}[H]
\centering
\begin{tabular}{|c|c|} \hline
$\mathfrak{g}$ & $\Delta$\\ 
\hline
$\Big. A(m-1,n-1)$
&$\delta_1-\delta_2, \dots, \delta_{n-1}-\delta_n, \delta_n-\epsilon_1, 
\epsilon_1-\epsilon_2, \dots, \epsilon_{m-1}-\epsilon_{m}$ \\
$\Big.B(m,n)$
&$\delta_1-\delta_2, \dots, \delta_{n-1}-\delta_n, \delta_n-\epsilon_1, 
\epsilon_1-\epsilon_2, \dots, \epsilon_{m-1}-\epsilon_m, \epsilon_m$ \\
$\Big. B(0,n)$
&$\delta_1-\delta_2, \dots, \delta_{n-1}-\delta_n, \delta_n$ \\
$\Big. C(n)$
&$\epsilon-\delta_1, \delta_1-\delta_2, \dots, \delta_{n-1}-\delta_n, 2\delta_n$ \\
$\Big. D(m,n)$
&$\delta_1-\delta_2, \dots, \delta_{n-1}-\delta_n, \delta_n-\epsilon_1, 
\epsilon_1-\epsilon_2, \dots, \epsilon_{m-1}-\epsilon_m, \epsilon_{m-1}+\epsilon_m$ \\ 
$\Big. D(2,1;\alpha)$
&$\epsilon_1-\epsilon_2-\epsilon_3, 2\epsilon_2, 2\epsilon_3$ \\
$\Big. F(4)$
&$\frac{1}{2} (\delta-\epsilon_1-\epsilon_2-\epsilon_3), \epsilon_3, \epsilon_2-\epsilon_3, 
\epsilon_1-\epsilon_2$ \\
$\Big. G(3)$
&$\delta+\epsilon_3, \epsilon_1, \epsilon_2-\epsilon_1$ \\

\hline
\end{tabular}
\end{table}

If $\Delta=\{\alpha_1,\dots,\alpha_r\}$ is a distinguished simple system of a basic classical Lie superalgebra of type I (resp. type II), then $\Delta_{\bar{0}}=\{\alpha_i:  i\neq s\}$ (resp. $\Delta_{\bar{0}}=\{\alpha_i, \gamma: i\neq s\}$) is a simple system for $\lie g_{\bar{0}}$, where $\gamma=\sum_{i=s}^{r}c_i\alpha_i$ with labels $c_i$ as in \cite[Table 1]{K78}. We have that $\gamma$ is the longest simple root of $C_n$ in the case of $B(m,n)$ and $D(m,n)$ and the simple root of $A_1$ in the case of $F(4), G(3)$ and $D(2, 1;\alpha)$. 
\subsection{}\label{section25}
Let $D=\mbox{diag}(d_i)_{i\in I}$ and $B=(b_{i,j})_{i,j\in I}$ be diagonal and symmetric matrices such that $A=DB$. We recall the notion of a Chevalley basis; for more details we refer the reader to \cite{FG11} and \cite{IK01}. By \cite[Theorem 3.9]{IK01} we can choose for any basic classical Lie superalgebra a homogeneous vector space basis $\{x_{\alpha}, h_i : i\in I, \alpha\in R\}$ consisting of root vectors $x_{\alpha}\in \lie g_{\alpha}, \alpha\in R$ such that the following holds:
\begin{equation*}\label{0} \big\{h_1,\dots,h_n\big\} \mbox{ is a basis of }\lie h\ \text{with } \alpha_i(h_j)=a_{j,i},
\end{equation*}
\begin{equation*}\label{1} [h_i,h_j]=0,\quad [h_i,x_{\alpha}]=\alpha(h_i)x_{\alpha},\quad [x_{\alpha},x_{-\alpha}]=\sigma_{\alpha}h_{\alpha},\quad \forall i,j\in I,\ \alpha\in R\end{equation*}
\begin{equation}\label{2}[x_{\alpha},x_{\beta}]=C_{\alpha,\beta}\ x_{\alpha+\beta}, \quad \forall \alpha,\beta\in R \ \mbox{ with $\alpha+\beta\in R$ and $C_{\alpha,\beta}\in \mathbb{Z}\backslash \{0\}$},
\end{equation}
where $\sigma_\alpha,h_\alpha$ are defined as follows. We have $\sigma_\alpha=-1$ if $\alpha\in R_{\bar{1}}^-$ and $\sigma_\alpha=1$ otherwise. For a root $\alpha=\sum_{i=1}^nk_i\alpha_i\in R$ we define its coroot 
$$h_{\alpha}= d_{\alpha}\sum_{i=1}^n k_id_i^{-1}h_i,\ \text{ where }\ d_{\alpha}=\begin{cases}
\frac{2}{(\alpha,\alpha)},& \text{if $(\alpha,\alpha)\neq 0$}\\
d_s,& \text{if $(\alpha,\alpha)=0$.}
\end{cases}$$
There are some further restrictions on the structure constants $C_{\alpha,\beta}$ which are not important in the remainder of this paper; see for example \cite[Definition 3.3]{FG11} or \cite[Theorem 3.9]{IK01}. 

\subsection{}
We recall the Poincaré-Birkoff-Witt theorem for Lie superalgebras. We denote by $\mathbf{T}(\mathfrak{g})$ the tensor superalgebra and let $\mathbf{J}$ the ideal of $\mathbf{T}(\mathfrak{g})$ generated by the elements of the form 
$$[x,y]-x\otimes y+(-1)^{|x||y|}y\otimes x,\ \ x,y\in\lie g,$$
where $|x|$ denotes the parity of $x$. The universal enveloping algebra $\mathbf{U}(\mathfrak{g})$ is defined as the quotient $\mathbf{T}(\mathfrak{g})/\mathbf{J}$. Note that the ideal $\mathbf{J}$ is graded and hence $\mathbf{U}(\mathfrak{g})$ is an associative superalgebra.
If $x_1,\dots,x_p$ is a vector space basis of $\mathfrak{g}_{\bar{0}}$ and $y_1,\dots,y_q$ a vector space basis of $\mathfrak{g}_{\bar{1}}$, then the PBW theorem says that the set of monomials
$$y_{i_1}\cdots y_{i_\ell} x_1^{k_1}\cdots x_p^{k_p},\ \ k_i\geq 0,\ \ 1\leq i_1<\cdots<i_\ell\leq q$$
forms a basis of $\mathbf{U}(\mathfrak{g})$. We define a filtration on $\mathbf{U}(\mathfrak{g})$:
\begin{equation}\label{filt1}\mathbf{U}(\mathfrak{g})_0\subseteq \mathbf{U}(\mathfrak{g})_1\subseteq \cdots\subseteq \mathbf{U}(\mathfrak{g})_p\subseteq \cdots\end{equation}
where $\mathbf{U}(\mathfrak{g})_0= \mathbb{C}\cdot 1$ and $\mathbf{U}(\mathfrak{g})_p$ is generated by the products of the form 
$$a_1\cdots a_m,\ \ 0\leq m\leq p,\ a_i\in \lie g.$$
We call \eqref{filt1} the PBW filtration of the universal enveloping algebra and the associated graded space with respect to \eqref{filt1}
$$\bigoplus_{k\in\bz} \mathbf{U}(\mathfrak{g})_k/\mathbf{U}(\mathfrak{g})_{k-1}$$
admits an obvious $\bz$-graded algebra structure obtained from that in $\mathbf{U}(\mathfrak{g})$ by going to the quotients. The structure of this algebra can be realized as follows. Define the symmetric superalgebra
$$\Lambda=S(\lie g_{\bar{0}})\otimes \Lambda(\lie g_{\bar{1}}),$$
to be the tensor product of the symmetric algebra of $\lie g_{\bar{0}}$ with the exterior algebra of $\lie g_{\bar{1}}$ with grading
$$\Lambda_n:=\bigoplus_{m=0}^n \left(S_m(\lie g_{\bar{0}})\otimes \Lambda^{(n-m)}(\lie g_{\bar{1}})\right).$$
The following is a straightforward consequence of the PBW theorem and can be found for example in \cite[Chapter 2]{Sch97}.
\begin{prop}
The associated graded space of $\mathbf{U}(\mathfrak{g})$ with respect to the PBW filtration \eqref{filt1} is isomorphic to the algebra $S(\lie g_{\bar{0}})\otimes \Lambda(\lie g_{\bar{1}})$ as $\bz$--graded superalgebras.

\hfill\qed
\end{prop}
\subsection{} Now we discuss finite-dimensional representations of basic classical Lie superalgebras and the PBW filtration. For $\lambda\in \lie h^{*}$ we define a one--dimensional irreducible $\lie b$--module $\mathbb{C}_{\lambda}:=\mathbb{C}v_\lambda$ by
$$\lie n^+v_{\lambda}=0,\ hv_{\lambda}=\lambda(h)v_{\lambda}, \ \forall h\in \lie h.$$
The module $L(\lambda):=\mathbf{U}(\lie g)\otimes_{\mathbf{U}(\lie b)}\mathbb{C}_{\lambda}$ contains a unique maximal submodule $J(\lambda)$. It is clear that the quotient $V(\lambda)=L(\lambda)/J(\lambda)$ is an irreducible representation with $V=\mathbf{U}(\lie n^-)v_{\lambda}$ (for simplicity we denote the highest weight vector $1\otimes v_{\lambda}$ also by $v_{\lambda}$). 
The following proposition is stated in \cite[Proposition 2.2]{K78}.
\begin{prop}Let $V$ be a finite-dimensional irreducible module for the basic classical Lie superalgebra $\lie g$. Then there exists $\lambda\in \lie h^*$ such that $V\cong V(\lambda)$. Moreover, $V(\lambda)\cong V(\mu)$ if and only if $\lambda=\mu$.
\hfill\qed
\end{prop}
The construction of finite-dimensional irreducible representations for basic classical Lie superalgebras is quite similar as for simple Lie algebras, but the parametrizing set again depends on the choice of the Borel subalgebra. Since we have assumed in this paper that $\lie b$ is distinguished we can describe the Zariski dense set $P^+=\{\lambda\in \lie h^{*}\mid \dim V(\lambda)<\infty\}$ as follows. Let $P_{\bar{0}}^+$ be the set of dominant integral weights for $\lie g_{\bar{0}}$. One of the necessary conditions for $V(\lambda)$ to be finite-dimensional is that $\lambda(h_\alpha)\in\bz_+$ for all $\alpha\in R_{\bar{0}}^+$, i.e. $P^+\subseteq P_{\bar{0}}^+$. For the special linear Lie superalgebra this condition is also sufficient, but there are a few extra conditions in the remaining types. All characterizing properties can be found in \cite[Proposition 2.3]{K78}. 
\subsection{}We give a more explicit construction of the finite-dimensional irreducible $\lie g$-modules in the typical case. Recall that $V(\lambda)$ is called typical if $(\lambda+\rho,\alpha)\neq 0$ for all $\alpha\in R_{\bar{1}}^+$. For a complete and explicit characterization of typical representations we refer to \cite[pg. 620-622]{K78}. The following proposition stated in \cite[Theorem 1]{K78} gives generators and relations for typical finite-dimensional irreducible $\lie g$--modules. Recall the definition of $\gamma$ from Section~\ref{section24}.
\begin{prop}\label{prop1}Let $\lambda\in P^+$ a typical weight.
We have an isomorphism of $\lie g$-modules
$$V(\lambda)\cong \mathbf{U}(\lie g)/M(\lambda),$$
where $M(\lambda)$ is the left ideal generated by $\lie n^+$, $(h-\lambda(h)\cdot 1)$ for all $h\in \lie h$ and
$$\begin{cases} (x_{-\alpha_i})^{\lambda(h_i)+1}\mbox{ for } i\neq s,& \text{ if $\lie g$ is of type I}\\
(x_{-\alpha_i})^{\lambda(h_i)+1} \mbox{ for } i\neq s,\ (x_{-\gamma})^{2\frac{(\lambda,\gamma)}{(\gamma,\gamma)}+1},& \text{ if $\lie g$ is of type II.}
\end{cases}
$$
Moreover,
\begin{equation}\label{dimf}
 \dim V(\lambda)= 2^{|R_{\bar{1}}^+|}\prod_{\alpha\in R_{\bar{0}}^+}\frac{(\lambda+\rho,\alpha)}{(\rho_0,\alpha)}.\end{equation}
 \hfill\qed
\end{prop}
\begin{rem}\label{grad} ${ }$
\begin{enumerate}
\item All finite-dimensional irreducible representations of $B(0,n)$ are typical.
\item There is a distinguished $\bz$-gradation $\lie g=\bigoplus_{i\in\bz} \lie g_i$ on $\lie g$ by putting
$$\text{deg}(h_i)=0, \text{deg}(e_i)=0=\text{deg}(f_i),\ \ i\neq s,\ \ \text{deg} (e_s)=-\text{deg} (f_s)=1.$$
Note that this grading is compatible with the $\bz_2$-grading, i.e. $\lie g_{\bar{0}}$ is the direct sum over all even homogeneous spaces. Furthermore, $\lie g_i=0$ for all $|i|>1$ (resp. $|i|>2$) if $\lie g$ is of type I (resp. type II). Set $\lie g_+=\bigoplus_{i\in\bz_+} \lie g_i$ and let $V_{\lie g_{0}}(\lambda)$ the finite--dimensional irreducible $\lie g_{0}$--module of highest weight $\lambda$. We  extend $V_{\lie g_{0}}(\lambda)$ to a $\mathbf{U}(\lie g_+)$-module by requiring that $\lie g_i$ acts as zero for all $i>0$. The Kac module is then defined as 
$$K(\lambda)\cong \text{Ind}^{\lie g}_{\lie g_+} V_{\lie g_{0}}(\lambda)/M,$$
where $M=0$ if $\lie g$ is of type I and $M=\mathbf{U}(\lie g)(x_{-\gamma})^{2\frac{(\lambda,\gamma)}{(\gamma,\gamma)}+1}v_{\lambda}$ otherwise. One of the fundamental results in the representation theory of basic classical Lie superalgebras is that $K(\lambda)$ is irreducible if and only if $\lambda$ is typical. 
\item For atypical representations we have in general that $\dim V(\lambda)$ is strictly less than the right hand side of \eqref{dimf}. The equality in \eqref{dimf} is essential for the rest of the paper.
\end{enumerate}
\end{rem}
\subsection{}Recall the PBW filtration from \eqref{filt1}. We define an induced filtration on $V(\lambda), \lambda\in P^+$ as follows:
\begin{equation}\label{filt2}\mathbf{U}(\mathfrak{n}^-)_0 v_{\lambda}\subseteq \mathbf{U}(\mathfrak{n}^-)_1v_{\lambda}\subseteq \cdots\subseteq \mathbf{U}(\mathfrak{n}^-)_p v_{\lambda}\subseteq \cdots\subseteq V(\lambda).
\end{equation}
The associated graded space with respect to \eqref{filt2} is defined as
$$\gr V(\lambda):=\bigoplus_{k\in\bz} V^k(\lambda)/V^{k-1}(\lambda),\ \ V^k(\lambda):=\mathbf{U}(\mathfrak{n}^-)_kv_{\lambda}$$
The following lemma is straightforward.
\begin{lem}The action of $\mathbf{U}(\lie n^-)$ on $V(\lambda)$ induces an action of $S(\lie n^-_{\bar{0}})\otimes \Lambda(\lie n^-_{\bar{1}})$ on $gr V(\lambda)$, which turns $\gr V(\lambda)$ into a cyclic representation. Moreover, there is an induced action of $\mathbf{U}(\lie n^+)$ on the associated graded space $\gr V(\lambda)$.
\hfill\qed
\end{lem}
Let $\mathbf{I}(\lambda)$ the left ideal of $S(\lie n^-_{\bar{0}})\otimes \Lambda(\lie n^-_{\bar{1}})$ such that 
\begin{equation}\label{viaid}
\gr V(\lambda)\cong S(\lie n^-_{\bar{0}})\otimes \Lambda(\lie n^-_{\bar{1}})/\mathbf{I}(\lambda).\end{equation}
We will make the action of $\mathbf{U}(\lie n^+)$ on $\gr V(\lambda)$ via the identification \eqref{viaid} more explicit. Define differential operators $\partial_{\alpha}$, $\alpha\in R^{+}$ on $S(\lie n^-_{\bar{0}})\otimes \Lambda(\lie n^-_{\bar{1}})$ by:
$$
\partial_{\alpha} x_{-\beta} :=
\begin{cases}
x_{-\beta + \alpha}, \, &\textrm{if} \, \beta - \alpha \in R^+\\
0, \,& \textrm{else.}
\end{cases}$$
Then $x_{\alpha}$, $\alpha\in R^+$ acts on the associated graded space via the differential operator $\partial_\alpha$. The goal of this paper is to make a first step towards understanding the structure of $\gr V(\lambda)$ provided that the structure of $\gr V_{\lie g_{\bar{0}}}(\lambda)$ is known, where $V_{\lie g_{\bar{0}}}(\lambda)$ is the finite-dimensional irreducible $\lie g_{\bar 0}$-module and $\gr V_{\lie g_{\bar{0}}}(\lambda)$ is defined as in \cite{FFL11}. The structure of $\gr V_{\lie g_{\bar{0}}}(\lambda)$ including a computation of a monomial basis parametrized by the lattice points of a convex polytope has been worked out in \cite{FFL11} for type $A_n$, in \cite{FFL12} for type $C_n$, in type $B_3$ in \cite{BK19} and for type $G_2$ in \cite{G15}.
\subsection{}As a first step we will reduce the problem of finding a PBW basis into several pieces. We denote by $\lie g(1),\dots \lie g(p)$ the simple finite--dimensional Lie algebras which appear in the semi--simple part of the reductive Lie algebra $\lie g_{\bar{0}}$ and let $R^+(i)\subset R^+_{\bar{0}}$ the corresponding set of positive roots. Let $\lie g$ a Lie superalgebra of type II. We will assume without loss of generality that $\gamma\in R^+(1)$ (e.g. if $\lie g$ is of type $B(m,n)$, then $\lie g(1)$ is of type $C_n$ and $\lie g(2)$ of type $B_m$). Since $P^+\subseteq P_{\bar{0}}^+$ we can write any $\lambda\in P^+$ as a linear combination \begin{equation}\label{32}\lambda=\lambda_1+\cdots+\lambda_{p},\end{equation}
where $\lambda_i$ is a dominant integral weight for the simple Lie algebra $\lie g(i)$. The following lemma is straightforward to check by using the fact that 
$$(\rho_{\bar{1}}, \alpha)=0,\ \ \forall \alpha\in R^+(i),\ i\neq 1$$
and gives a factorization of the dimension formula \eqref{dimf}.
\begin{lem}\label{dimzer}Let $\lambda\in P^+$ a typical weight and $\lambda=\lambda_1+\cdots+\lambda_p$ a decomposition as in \eqref{32}. Then
 $$\dim V(\lambda)= \dim V(\lambda_1) \left(\prod_{i=2}^{p} \dim V_{\lie g(i)}(\lambda_i)\right).$$
 \hfill\qed
\end{lem}
We will need some more notation to state our reduction result for type II Lie superalgebras (see Proposition~\ref{reduction}). Let $\lie n^-(i)$ be the Lie superalgebra generated by the root vectors $x_{-\alpha},\ \alpha\in R^+(i)$. Since $v_{\lambda}$ satisfies the relations of $V_{\lie g(i)}(\lambda_i)$ by Proposition~\ref{prop1}, we immediatly obtain 
\begin{equation}\label{red0}V_{\lie g(i)}(\lambda_i)\cong \mathbf{U}(\lie n^-(i))v_{\lambda}\subseteq V(\lambda).\end{equation}
Moreover, it is not hard to check that $(\lie n^-(1)\oplus \lie n^-_{\bar{1}})$ is a Lie superalgebra. We also get 
\begin{equation}\label{red01}V_{}(\lambda_1)\cong \mathbf{U}(\lie n^-(1)\oplus \lie n^-_{\bar{1}})v_{\lambda}\subseteq V(\lambda).\end{equation}
\begin{prop}\label{reduction}Let $\mathcal{B}_i$ a PBW basis for the space $\gr V_{\lie g(i)}(\lambda_i)$ for $2\leq i \leq n$ and $\mathcal{B}_1$ a basis for the space $\gr V(\lambda_1)$. Then, we have that 
\begin{equation}\label{red1}\{b_1\cdots b_p v_{\lambda}: b_i\in \mathcal{B}_i\}\end{equation}
forms a basis for $\gr V(\lambda)$. 
\proof
The cardinality of \eqref{red1} coincides with $\dim V(\lambda)$ by Lemma~\ref{dimzer}. It remains to show that the above set is a generating set. Clearly, the set of elements 
$$u_1\cdots u_pv_{\lambda},\ \ \text{where }  u_i\in \mathbf{U}(\lie n^-(i)) \text{ is a monomial}$$
spans $\gr V(\lambda)$. From \eqref{red0} and \eqref{red01} we obtain that 
$$\gr V_{\lie g(i)}(\lambda_i)\cong S(\lie n^-(i))v_{\lambda},\ (2\leq i\leq n),\ \ \gr V(\lambda_1)\cong\left( S(\lie n^-(1))\otimes \Lambda( \lie n^-_{\bar{1}})\right)v_{\lambda}.$$ 
Now we can rewrite $u_pv_{\lambda}$ as a linear combination of elements $b_pv_{\lambda},\ b_p\in \mathcal{B}_p$. By the commutativity we can pass $u_{p-1}$ through all elements in $\mathcal{B}_p$ and write $u_{p-1}v_{\lambda}$ in the basis $\mathcal{B}_{p-1}$. By repeating the above steps the claim follows.
\endproof
\end{prop}
We emphasize the importance of Proposition~\ref{reduction}. In order to determine a PBW basis for $\gr V(\lambda)$ we only have to compute a PBW basis for the representations of the underlying simple Lie algebras (most of the cases are known in the literature; see for example \cite{BK19,FFL11,FFL12,G15}) and a PBW basis for $\gr V(\lambda_1)$, where $\lambda_1$ is a dominant integral weight for the simple Lie algebra $\lie g(1)$. With other words, we can lift a PBW basis of $\gr V_{\lie g_{\bar{0}}}(\lambda)$ 
to a PBW basis of $\gr V(\lambda)$ provided that the structure of $\gr V(\lambda_1)$ is known. An application of this important result can be found for example in the Appendix, when we contruct a basis for the exceptional cases. 
\section{Lifting PBW bases: type I case}\label{section3}
 The  answer to the natural question whether we can lift a PBW basis of $\gr V_{\lie g_{\bar{0}}}(\lambda)$ to a PBW basis of  $\gr V(\lambda)$ for type I basic classical Lie superalgebras turns out to be quite easy.
\subsection{}We enumerate the positive odd roots $R_{\bar{1}}^+=\{\beta_1,\dots,\beta_d\}$. The following theorem gives a PBW basis for the associated graded space.
\begin{thm}\label{thm-typeI} Let $\lie g$ be a basic classical Lie superalgebra of type I and $\lambda\in P^+$ a typical weight. Further, let $\mathcal{B}=\{b_1,\dots, b_\ell\}$ a PBW basis of $\gr V_{\lie g_{\bar{0}}}(\lambda)$. Then the set
$$\mathcal{B}=\{x^{\mathbf{r}} b_i v_{\lambda}: \ \ 1\leq i\leq \ell,\ \ \mathbf{r}\in \{0,1\}^{d}\}$$
forms a PBW basis of $\gr V(\lambda)$, where $x^{\mathbf{r}}:=x_{-\beta_1}^{r_1}\cdots x_{-\beta_d}^{r_d}$.
\end{thm}
\proof Since $V(\lambda)\cong \mathbf{U}(\lie g_{-1})\otimes V_{\lie g_{\bar{0}}}(\lambda)$ as a vector space by Remark \ref{grad}, we have the correct cardinality of $\mathcal{B}$. Hence it will be enough to show that the elements of $\mathcal{B}$ are linearly independent. Assume by contradiction that we have the linearly dependence in $\gr V(\lambda)$: 
\begin{equation}\label{dep}\sum_{i, \mathbf{r}}\lambda_{i,\mathbf{r}} \ x^{\mathbf{r}}b_i v_{\lambda}=0,\ \ \lambda_{i,\mathbf{r}}\in\bc, \end{equation}
where we can assume without loss of generality that all summands have the same PBW degree, i.e. the linear dependence is in the space $V^k(\lambda)/V^{k-1}(\lambda)$ for some $k\in\bn$. This means that the left hand side of \eqref{dep} is contained in $ V^{k-1}(\lambda),$ say it equals an element of the form $Z=\sum_{j}\mu_{j}  z_jv_{\lambda}$ with $z_j\in \mathbf{U}(\lie n^-)_{k-1}, \mu_j\in\bc$. Writing each $z_j$ in PBW order (i.e. as a product of elements in $\lie n^{-}_{\bar{1}}$ followed by a product of elements in $\lie n^{-}_{\bar{0}}$) we can assume that $Z$ is of the following form 
$$Z=\sum_{\mathbf{p}}\mu_{\mathbf{p}} \ x^{\mathbf{p}}v_{\mathbf{p}} v_{\lambda},\ \ \mu_{\mathbf{p}}\in\bc$$
for some elements $v_{\mathbf{p}}\in \mathbf{U}(\lie n_{\bar{0}}^-)$. Now using the fact that  
$$V(\lambda)\cong \mathbf{U}(\lie g)\otimes_{\mathbf{U}(\lie g_+)} V_{\lie g_{\bar{0}}}(\lambda)\cong \mathbf{U}(\lie g_{-1})\otimes_{\bc} V_{\lie g_{\bar{0}}}(\lambda)$$ we get
$$0=\sum_{i, \mathbf{r}}\lambda_{i,\mathbf{r}} \ x^{\mathbf{r}}\otimes b_i v_{\lambda}- \sum_{\mathbf{p}}\mu_{\mathbf{p}} \ x^{\mathbf{p}}\otimes v_{\mathbf{p}}v_{\lambda}.$$ Since $\{x^{\mathbf{r}}\}_{\mathbf{r}\in\{0,1\}^d}$ is a linearly independent subset of $\mathbf{U}(\lie g_{-1})$ we obtain by collecting all the coefficients of $x^{\mathbf{r}}$ that
$\lambda_{1,\mathbf{r}}b_1+\cdots+\lambda_{\ell,\mathbf{r}}b_{\ell}=\mu_{\mathbf{r}}v_{\mathbf{r}},$
which means that $\lambda_{1,\mathbf{r}}b_1+\cdots+\lambda_{\ell,\mathbf{r}}b_{\ell}=0$ in the space $\gr V_{\lie g_{\bar{0}}}(\lambda)$. This is a contradiction and the claim follows.
\endproof
Recall from \eqref{viaid} the definition of the left ideal $\mathbf{I}(\lambda)$.
\begin{cor}\label{cor123}
We have that $\mathbf{I}(\lambda)$ is the left ideal in $S(\lie n^-_{\bar{0}})\otimes \Lambda(\lie n^-_{\bar{1}})$ generated by the elements 
\begin{equation}\label{idde2}\mathbf{U}(\lie n_{\bar 0}^+) \circ x^{2\frac{(\lambda,\alpha)}{(\alpha,\alpha)}+1}_{-\alpha}, \ \alpha\in R_{\bar{0}}^+.\end{equation}
Moreover, for all $\lambda,\mu\in P^+$ typical weights we have
$$\gr V(\lambda+\mu)\cong S(\lie n^-_{\bar{0}})\otimes \Lambda(\lie n^-_{\bar{1}})(v_{\lambda}\otimes v_{\mu})\subseteq \gr V(\lambda)\otimes \gr V(\mu).$$
\proof
The first part of the corollary follows by combining \cite[Theorem A]{FFL11} with Theorem~\ref{thm-typeI}. We denote by $S_{\lie g_{\bar{0}}}(\lambda)$ the lattice points of the convex polytope constructed in \cite[Definition 2]{FFL11} and \cite[Section 2]{FFL12} respectively. Recall the important Minkowski sum property of the polytope: $S_{\lie g_{\bar{0}}}(\lambda+\mu)=S_{\lie g_{\bar{0}}}(\lambda)+S_{\lie g_{\bar{0}}}(\mu)$. By \eqref{idde2} we immediately obtain a surjective map 
\begin{equation}\label{1223}\gr V(\lambda+\mu)\twoheadrightarrow S(\lie n^-_{\bar{0}})\otimes \Lambda(\lie n^-_{\bar{1}})(v_{\lambda}\otimes v_{\mu})\end{equation}
and the injectivity would follow from the fact that the set 
$$\{x^{\mathbf{r}} x^{\mathbf{s}}v_{\lambda}: \ \ {\mathbf{s}}\in S_{\lie g_{\bar{0}}}(\lambda+\mu),\ \ \mathbf{r}\in \{0,1\}^{d}\}$$
is linearly independent in the right hand side of \eqref{1223}. But this follows the same strategy as \cite[Proposition 6]{FFL11}. The key ingredient is the Minkowski property. 
\endproof
\end{cor}
\begin{rem}Another way to construct a basis for $V(\lambda)$ would be to compute the $\lie g_{\bar 0}$ decomposition and to take a basis for each direct summand. We emphasize that this basis will generically not be compatible with the PBW filtration, i.e. the image of this basis in $\gr V(\lambda)$ will not be a basis. For example, let $\mathfrak{g}=\mathfrak{sl}(3,1)$ and consider the typical representation $\lambda=\epsilon_1$. Then $x_{-(\epsilon_1-\delta_1)}x_{-(\epsilon_3-\delta_1)}v_{\lambda}$ will be a $\lie g_{\bar 0}$ highest weight vector (up to a filtration) of weight $\epsilon_1+\epsilon_2$. Hence $x_{-(\epsilon_2-\epsilon_3)}x_{-(\epsilon_1-\delta_1)}x_{-(\epsilon_3-\delta_1)}v_{\lambda}$ will be a basis vector of $V(\lambda)$, but it vanishes in the associated graded space since
$$x_{-(\epsilon_2-\epsilon_3)}x_{-(\epsilon_1-\delta_1)}x_{-(\epsilon_3-\delta_1)}v_{\lambda}=x_{-(\epsilon_1-\delta_1)}x_{-(\epsilon_2-\delta_1)}v_{\lambda}.$$ For another example consider the Lie superalgebra $\lie g=\mathfrak{osp}(1,4)$ and $\lambda=\delta_1+\delta_2$. Then $x_{-\delta_2}v_\lambda$ is a highest weight vector of highest weight $\delta_1$ and hence $x_{-(\delta_1-\delta_2)}x_{-\delta_2}v_{\lambda}$ would be a basis vector of $V(\lambda)$. But since $x_{-(\delta_1-\delta_2)}x_{-\delta_2}v_{\lambda}=x_{-\delta_1}v_{\lambda}$ we see that this vector vanishes in the associated graded space.
\end{rem}

\section{Lifting PBW bases: the $\mathfrak{osp}(1,2n)$ case}\label{section5}
Here we consider the Lie superalgebra $\mathfrak{osp}(1,2n)$. Recall that the set of simple roots for $\lie g_{\bar{0}}$ is given by $\{\alpha_1,\dots,\alpha_n\}$, where $\alpha_i=\delta_i-\delta_{i+1}, 1\leq i<n$ and $\alpha_n=2\delta_n$. All positive roots of $\lie g_{\bar{0}}$ are given by
\begin{align*}
\alpha_{i,j}=\alpha_i+\alpha_{i+1}+\dots +\alpha_j,\ 1\le i\le j\le n,\\
\alpha_{i, \overline{j}} = \alpha_i + \alpha_{i+1} + \ldots +
\alpha_n + \alpha_{n-1} + \ldots + \alpha_j, \ 1\le i\le j\le n.
\end{align*}
Recall also the set of positive odd roots $R_{\bar{1}}^+=\{\delta_1,\dots,\delta_n\}$ and let $\{\varpi_1,\dots,\varpi_n\}$ the set of fundamental weights for $\lie g_{\bar{0}}$.

\subsection{}We introduce the notion of an orthosymplectic Dyck path by slightly modifying the definition of a symplectic Dyck path from \cite[Definition 1.2]{FFL12}.
\begin{defn}\label{dyckpath}\rm
An orthosymplectic Dyck path is a sequence
\[
\mathbf{p}=(p(0), p(1),\dots, p(k)), \ k\ge 0
\]
of positive roots which satisfies one of the following conditions:
\begin{itemize}
\item[{\it (a)}] The sequence $\mathbf{p}$ is a symplectic Dyck path satisfying $p(0)=\alpha_i$ and $p(k)=\alpha_j$ for some $ 1 \leq i \leq j< n$;
\item[{\it (b)}] The sequence $(p(0), p(1),\dots, p(k-1))$ is a symplectic Dyck path satisfying $p(0)=\alpha_i$ and $p(k-1) = 2\delta_j$ for some $ i \le j \leq n$ and $p(k)=\delta_i$.
\end{itemize}
\end{defn}
Denote by $\mathcal{D}$ the set of all orthosymplectic Dyck paths. For a dominant integral $\lie g_{\bar{0}}$ weight $\lambda=\sum_{i=1}^n m_i \varpi_i$ and $\mathbf{p}\in\mathcal{D}$ we set 
$$M_{\mathbf{p}}(\lambda)=\begin{cases}m_i+\cdots +m_j, & \text{ if $\mathbf{p}$ satisfies (a)}\\m_i+\cdots +m_n,&\text{ otherwise. }
\end{cases}$$
Denote by $P(\lambda)\subseteq \mathbb{R}^{n(n+1)}_{\geq 0}$ the polytope
\begin{equation}
\label{polytopeequation}
P(\lambda):=\Big\{(s_{\alpha})_{\alpha\in R^+}\mid \forall \mathbf{p}\in \mathcal{D}:
s_{p(0)} + \dots + s_{p(k)}\le M_{\mathbf{p}}(\lambda),\ \forall \beta\in R_{\bar{1}}^+: s_{\beta}\leq 1 
\Big\},
\end{equation}
and let $S(\lambda)$ be the set of integral points in $P(\lambda)$. The main theorem of this section is the following.

\begin{thm}\label{thm1}Let $\lambda=\sum_{i=1}^n m_i\varpi_i\in P^+$. 
\begin{enumerate}
\item The monomials $x^{\mathbf{s}},\ \mathbf{s}\in S(\lambda)$
form a basis of $\gr V(\lambda)$ where $$x^{\mathbf{s}}:=\prod_{\alpha\in R^+}x_{-\alpha}^{s_\alpha}.$$
\item We have that $\mathbf{I}(\lambda)$ is the left ideal in $S(\lie n^-_{\bar{0}})\otimes \Lambda(\lie n^-_{\bar{1}})$ generated by the elements 
\begin{equation}\label{idde}\mathbf{U}(\lie n_{\bar 0}^+) \circ x^{m_i+\cdots+m_j+1}_{-(\delta_i-\delta_j)}, \ 1\leq i<j \leq n, \ \ \mathbf{U}(\lie n_{}^+) \circ x^{m_i+\cdots+m_n+1}_{-2\delta_i}, \ 1\leq i\leq n.\end{equation}
\item For $\lambda,\mu\in P^+$ we have
$$\gr V(\lambda+\mu)\cong S(\lie n^-_{\bar{0}})\otimes \Lambda(\lie n^-_{\bar{1}})(v_{\lambda}\otimes v_{\mu})\subseteq \gr V(\lambda)\otimes \gr V(\mu).$$
\end{enumerate}
\end{thm}
The proof of the above theorem will be given in the rest of this section.

\begin{rem}
Note that the elements $x^{\mathbf{s}}$ in Theorem~\ref{thm1} are only well-defined up to a sign. We could avoid this by choosing a total order on the set of odd roots and order the odd root vectors in $x^{\mathbf{s}}$ with respect to this order. For simplicity we will ignore this and write $x^{\mathbf{s}}$ without further comment.
\end{rem} 
\subsection{} In order to show that the lattice points of \eqref{polytopeequation} parametrize a basis of $\gr V(\lambda)$, we will need a staightening law similar to the one given in \cite[Theorem 2.4]{FFL12}. 
Recall the total order $\succ$ on the set of monomials in $S(\lie n^-_{\bar{0}})$ from \cite[Definition 2.6]{FFL12}. Similarly we can define a total order $>$ on the non-zero monomials in $\Lambda(\lie n_{\bar{1}}^-)$. Let
$$x_{-\delta_1}>x_{-\delta_2}>\cdots>x_{-\delta_n}$$ and let $>$ the induced homogeneous lexicographic ordering (which we also denote by $>$ for simplicity). We extend both to a total order on the set of monomials in $S(\lie n^-_{\bar{0}})\otimes \Lambda(\lie n^-_{\bar{1}})$ as follows. We define for a multi-exponent $\mathbf{s}\in \mathbb{Z}_+^{n^2}\times \{0,1\}^n$ the following elements
$$x^{\mathbf{s}}:=\prod_{\alpha\in R^+}x_{-\alpha}^{s_\alpha},\ \ x^{\mathbf{s}_{\bar{0}}}:=\prod_{\alpha\in R_{\bar{0}}^+}x_{-\alpha}^{s_\alpha},\ \ x^{\mathbf{s}_{\bar{1}}}:=\prod_{\alpha\in R_{\bar{1}}^+}x_{-\alpha}^{s_\alpha}.$$
We say $x^{\mathbf{s}}\succ x^{\mathbf{t}}$ if either 
\begin{itemize}
\item the total degree of $x^{\mathbf{s}}$ is greater than the total degree of $x^{\mathbf{t}}$
\item both have the same total degree, but $x^{\mathbf{s}_{\bar{1}}}> x^{\mathbf{t}_{\bar{1}}}$.
\item both have the same total degree, $s_\alpha=t_\alpha$ for all $\alpha\in R_{\bar{1}}^+$, but $x^{\mathbf{s}_{\bar{0}}}\succ x^{\mathbf{t}_{\bar{0}}}$. 
\end{itemize}
We will need the following lemma later.
\begin{lem}\label{strl}
\begin{enumerate}
\item Let $\mathbf{s}$ a multi-exponent such that $s_{-\delta_i}=0$ for all $1\leq i\leq n$. Moreover, assume that $\mathbf{s}$ is only supported on a symplectic Dyck path $\mathbf{p}$. Then there exists a homogeneous expression of the form 
$$x^{\mathbf s} + \sum_{\mathbf{t} \prec\, \mathbf{s}} c_{\mathbf{t}}  x^{\mathbf{t}} \in \mathbf{U}(\lie n_{\bar{0}}^+) \circ x^{|\mathbf{s}|}_{-\tau},\ \ c_{\mathbf{t}}\in\mathbb{C},$$
where $\tau=2\delta_i$ if $\mathbf{p}$ is a Dyck path satisfying $(b)$ and $\tau=\delta_i-\delta_j$ otherwise.
\item For every $k\in\mathbb{N}$ and $1\leq i\leq n$ we have the inclusion 
\begin{equation*}\left(\mathbf{U}(\lie n_{\bar 0}^+) \circ x^{k}_{-2\delta_i}\right)x_{-\delta_i}\subseteq \mathbf{U}(\lie n_{}^+) \circ x^{k+1}_{-2\delta_i}+ \sum_{p=i+1}^n \left(\mathbf{U}(\lie n_{\bar 0}^+) \circ x^{k}_{-2\delta_i}\right)x^{}_{-\delta_p}.\end{equation*}
\end{enumerate}
\begin{proof}
The first part of the lemma has been proved in \cite[Theorem 2.4]{FFL12}. In order to prove the second part let $\beta_1,\dots,\beta_r$ a sequence of positive even roots. Since $\delta_i-\alpha\notin R^+$ for all $\alpha\in R_{\bar{0}}^+$ unless $\alpha=\delta_i-\delta_{p}$ for some $p>i$ we obtain 
$$\partial_{\beta_1}\cdots \partial_{\beta_r}(x^{k}_{-2\delta_i}x^{}_{-\delta_i})-\left(\partial_{\beta_1}\cdots \partial_{\beta_r}(x^{k}_{-2\delta_i})\right)x^{}_{-\delta_i}\in \sum_{p=i+1}^n \left(\mathbf{U}(\lie n_{\bar 0}^+) \circ x^{k}_{-2\delta_i}\right)x^{}_{-\delta_p}.$$
Since $x^{k}_{-2\delta_i}x^{}_{-\delta_i}=\partial_{\delta_i}x^{k+1}_{-2\delta_i}$ the proof of the lemma is finished.
\end{proof}
\end{lem}
We claim that the first part of the above lemma holds even without the additional assumption $s_{-\delta_i}=0$. So let $\mathbf{s}$ a multi-exponent such that $s_{-\delta_i}\neq 0$ for some $i\in\{1,\dots,n\}$ and assume that $\mathbf{s}$ is supported on a Dyck path $\mathbf{p}$ with $p(0)=\alpha_i$ and $p(k)=2\delta_j$.
If we set $\mathbf{s}'$ as the multi--exponent obtained from $\mathbf{s}$ by setting $s_{-\delta_i}=0$ we obtain with Lemma~\ref{strl} the existence of a homogeneous expression
\begin{equation}\label{jjj}x^{\mathbf s'} + \sum_{\mathbf{t} \prec\, \mathbf{s}'} c_{\mathbf{t}}  x^{\mathbf{t}} \in \mathbf{U}(\lie n_{\bar 0}^+) \circ x^{|\mathbf{s}'|}_{-2\delta_i}.\end{equation}
If we multiply \eqref{jjj} by $x_{-\delta_i}$ and use Lemma~\ref{strl}(2) we obtain 
$$x^{\mathbf s} + \sum_{\mathbf{t} \prec\, \mathbf{s}} c_{\mathbf{t}}  x^{\mathbf{t}} \in \mathbf{U}(\lie n_{}^+) \circ x^{|\mathbf{s}|}_{-2\delta_i}+\sum_{p=i+1}^n \left(\mathbf{U}(\lie n_{\bar 0}^+) \circ x^{|\mathbf{s}'|}_{-2\delta_i}\right)x^{}_{-\delta_p}.$$
But the order on the monomials implies that each summand of $\left(\mathbf{U}(\lie n_{\bar 0}^+) \circ x^{k}_{-2\delta_i}\right)x^{}_{-\delta_p}$ is less than $x^{\mathbf{s}}$ and we obtain the desired straightening law which we summarize in the following corollary. 
\begin{cor}There exists a total order $\succ$ on the monomials in $S(\lie n^-_{\bar{0}})\otimes \Lambda(\lie n^-_{\bar{1}})$ such that
for any $\mathbf{s}\not\in S(\lambda)$ there exists a homogeneous expression
of the form
\begin{equation}\label{straighteninglaw}
x^{\mathbf s} -\sum_{\mathbf{t}\prec \mathbf s} c_{\mathbf t} x^{\mathbf t}
\end{equation}
which is contained in the left ideal generated by the elements \eqref{idde}.
\begin{proof}
If $\mathbf{s}$ is supported on a Dyck path, the statement is clear by the above discussion. Otherwise we use the fact that $\prec$ is a monomial order.
\end{proof}
\end{cor}
Now the same strategy as in the papers \cite{BK19,FFL11,FFL12} shows that the set $\{x^{\mathbf{s}}: \mathbf{s}\in S(\lambda)\}$ spans $S(\lie n^-_{\bar{0}})\otimes \Lambda(\lie n^-_{\bar{1}})/\mathbf{I}(\lambda)$ and hence the quotient space $\gr V(\lambda)$; we omit the details. In order to finish the proof of Theorem~\ref{thm1} we are left to show that $\{x^{\mathbf{s}}: \mathbf{s}\in S(\lambda)\}$ is a linearly independent subset of $\gr V(\lambda)$ (this shows part (1) and part (2) of the theorem) and that $\{x^{\mathbf{s}} (v_{\lambda}\otimes v_{\mu}): \mathbf{s}\in S(\lambda+\mu)\}$ is a linearly independent subset of $S(\lie n^-_{\bar{0}})\otimes \Lambda(\lie n^-_{\bar{1}})(v_{\lambda}\otimes v_{\mu})$ (this shows part (3)). This would follow from the following proposition; see for example the results of \cite{FFL17}.
\begin{prop}\label{hilfslem}Let $\lambda=\sum_{i=1}^n m_i\varpi_i\in P^+$.  
\begin{enumerate}
\item Let $i$ maximal with the property that $m_i\neq 0$. Then we have 
$$S(\lambda)\subseteq S(\lambda-\varpi_i)+S(\varpi_i).$$
In particular, for $\lambda,\mu\in P^+$ we have an equality
$$S(\lambda+\mu)=\big(S(\lambda)+S(\mu)\big)\cap \big(\mathbb{Z}_+^{n^2}\times \{0,1\}^n\big).$$
\item For each $i\in \{1,\dots,n\}$ the cardinality of $S(\varpi_i)$ is equal to $dim V(\varpi_i)$.
\end{enumerate}
\end{prop}
The proof of the above proposition can be found in the following two subsections. Given $\lambda = \sum_{i=1}^n m_i \varpi_i $ and $i$ maximal with the property that $m_i \neq 0$, we shall provide a procedure how to decompose a lattice point in the polytope $S(\lambda)$ as a sum of two elements in $S(\lambda-\varpi_i)$ and $S(\varpi_i)$ respectively. The procedure will depend on a total order. The equality in part (1) of the proposition is then an immediate consequence.
\subsection{Proof of Proposition~\ref{hilfslem}(1)}Define the total order on the set of positive roots
\begin{align}\label{ord1}
\alpha_{n,n} > \delta_n > \alpha_{n-1, \overline{n-1}} > \alpha_{n-1, n} > \alpha_{n-1, n-1} >& \delta_{n-1} > \alpha_{n-2, \overline{n-2}} > \cdots &\\&\notag\cdots >  \alpha_{n-2, n-2} > \delta_{n-2} > \ldots  > \alpha_{1,1} > \delta_1\end{align}
and consider the induced lexicographic ordering on the multi--exponents. Given a multi--exponent $\mathbf{s}\in S(\lambda)$ we define 
$$\mathbf{s}^1 := \max \{ \mathbf{t} \in S(\varpi_i) \mid t_\alpha \leq s_{\alpha}, \text{ for all } \alpha \in R^+\},
$$
where the maximum is with respect to the total order \eqref{ord1}. We claim that $\mathbf{s}-\mathbf{s}^1\in S(\lambda-\varpi_i)$ which would finish the proof of the first part of the proposition. Let $\mathbf{p}$ be a orthosymplectic Dyck path and assume first that $\mathbf{p}$ starts and ends at a simple even roots, $\alpha_s$ and $\alpha_e$ with $s\leq e$, different from $2 \delta_n$. We have to show that
\begin{equation}\label{11}
\sum_{\beta \in \mathbf{p}} (s_{\beta}-s^1_{\beta}) \leq m_s + \ldots + m_e-\delta_{i,\{s,\dots,e\}},
\end{equation}
where $\delta_{i,\{s,\dots,e\}}=1$ if $i\in\{s,\dots,e\}$ and $0$ otherwise. If $e<i$, there is nothing to be checked as the right hand side won't be changed by turning from $\lambda$ to $\lambda-\varpi_i$. Since $i$ is maximal with $m_i\neq 0$ we can assume that $e = i$. Let $\sigma$ be the maximal root (with respect to \eqref{ord1}) such that  $\sigma$ appears in the sequence $\mathbf{p}$ and $s_{\sigma} \neq 0$. If $s^1_{\sigma}\neq 0$, then we immediately get \eqref{11}; so assume without loss of generality that $s^1_{\sigma}=0$. By the maximality of $\mathbf{s}^1$, this is only possible if there exists a root $\tau>\sigma$ such that $s^1_{\tau}\neq 0$ and the multi-exponent $\mathbf{s}''$ with $s''_{\tau}=s''_{\sigma}=1$ and $s''_{\alpha}=0$ for $\alpha\notin\{\sigma,\tau\}$ is not contained in $S(\varpi_i)$. With other words, there is a Dyck path which contains both roots $\tau$ and $\sigma$. So we can replace the path $\mathbf{p}$ by a Dyck path $\mathbf{p'}$ which ends at $p'(k)=\delta_s$ (recall Definition~\ref{dyckpath}(b)) with the property $\tau\in\mathbf{p'}$ and $\beta\in \mathbf{p'}$ for all $\beta\in \mathbf{p}$ with $s_\beta\neq 0$. Hence
\[
\sum_{\beta \in \mathbf{p}} s_{\beta} < \sum_{\beta \in \mathbf{p}} s_{\beta}+s_{\tau}\leq \sum_{\beta \in \mathbf{p}'} s_\beta \leq m_s + \ldots + m_{i-1}+m_i 
\]
and the claim follows. Similarly, we can treat paths ending in $2\delta_j$. As before, we consider the maximal element on the path which is non-zero. Then either the corresponding element in $\mathbf{s}^1$ is also non-zero or there is a non-zero element in $\mathbf{s}^1$ which corresponds to a larger element than this maximal element.  
\begin{rem} Let $\mathbf{s}\in S(\lambda)$ and let $\tilde{\mathbf{s}}\in S(\lambda)$ the multi-exponent defined by 
$$\tilde{s}_{\alpha}=\begin{cases}s_{\alpha},&\text{if $\alpha\in R_{\bar{0}}^+$}\\
0,& \text{else}.
\end{cases}$$
Hence by \cite[Lemma 4.5]{FFL12} (the same result holds for $i$ maximal) there exists $\tilde{\mathbf{s}}^1\in S(\varpi_i)$ with $\tilde{s}^1_{\alpha}=0$ for all $\alpha\in R_{\bar{1}}^+$ such that $\tilde{\mathbf{s}}-\tilde{\mathbf{s}}^1\in S(\lambda-\varpi_i)$. Our desired element $\mathbf{s}^1\in S(\varpi_i)$ from Proposition~\ref{hilfslem}(1) is then given by
$$s^1_{\alpha}=\begin{cases}\tilde{s}^1_{\alpha},&\text{if $\alpha\in R_{\bar{0}}^+$}\\
1,& \text{if $s_{\alpha}\neq 0$ and $\alpha\in \{\delta_p: p\geq k+1\}$}\\
0,& \text{else}
\end{cases}$$
where $k\in\{1,\dots,n\}$ is the maximal column (in the Hasse diagram) such that there is a non-zero entry in $\tilde{\mathbf{s}}^1$.
\end{rem}
\subsection{Proof of Proposition~\ref{hilfslem}(2)}Here we prove that the cardinality of $S(\varpi_i)$ coincides with $\dim V(\varpi_i)$. For this we will use the notion of KT-tableaux; see for example \cite[Section 4]{KT90}. Let 
$J=\{\#_1,1,\bar{1},\#_2,2,\bar{2},\dots,\#_n,n,\bar{n}\}$ with ordering
$$\#_1<1<\bar{1}<\#_2<2<\bar{2}<\dots<\#_n<n<\bar{n}.$$
A KT-tableaux of shape $\epsilon\varpi_i$ ($\epsilon=2$ if $i=n$ and $\epsilon=1$ otherwise) is a column with $i$ boxes filled with entries from $J$ such that the entries from top to bottom are strictly increasing and the entry in row $r$ is greater or equal to $\#_r$. We denote the set of KT-tableaux of shape $\epsilon\varpi_i$ by $KT(\epsilon\varpi_i)$. For $\mathbf{s}\in S(\varpi_i)$ we define a KT-tableaux $T(\mathbf{s})\in KT(\epsilon\varpi_i)$ as follows. Since $\mathbf{s}\in S(\varpi_i)$ we must have the following: if $s_{\alpha}\neq 0$, then $\alpha=\alpha_{r,s}$ or $\alpha=\alpha_{r,\bar{s}}$ for some $r\leq i$. We take the tableaux 
$${\scriptsize\young(1,2,\vdots,i)}$$
and replace $r$ by $s+1$ if $s_{\alpha_{r,s}}\neq 0$ (we understand $n+1=\bar{n}$) and $r$ by $\bar{s}$ if $s_{\alpha_{r,\bar{s}}}\neq 0$. We denote the resulting tableaux (after permuting the entries in increasing order) by $Q(\mathbf{s})$. If $s_{\delta_j}\neq 0$, it is clear that $j\leq i$ and $s_{\alpha_{r,\bullet}}=0$ for all $r\geq j$. Hence there exists an entry $j$ in $Q(\mathbf{s})$. We replace $j$ by $\#_j$ in $Q(\mathbf{s})$ and define the resulting tableaux as $T(\mathbf{s})$. We claim that $T(\mathbf{s})$ is a KT-tableaux, where it remains to check that the entry in row $r$ is greater or equal to $\#_r$. For that it will be enough to prove that the entry in row $j$ of $Q(\mathbf{s})$ is contained in the set $\{j,\bar{j},\dots,n,\bar{n}\}$ for all $1\leq j\leq i$. But this is a straightforward calculation. 
\begin{lem}\label{inj}The map
$$S(\omega_i)\rightarrow KT(\epsilon\omega_i),\ \ \mathbf{s}\mapsto T(\mathbf{s}),$$
is injective.
\begin{proof}
Assume that $\mathbf{s}\neq \mathbf{s}'$ and $T(\mathbf{s})=T(\mathbf{s}')$. It follows immediately
\begin{equation}\label{1234}|\{1\leq r\leq n: s_{\delta_r}=1\}|=|\{1\leq r\leq n: s'_{\delta_r}=1\}|.\end{equation}
We first suppose that there exists $j\in\{1,\dots,n\}$ such that $s_{\delta_j}=1$ but $s'_{\delta_j}=0$ and assume that $j$ is minimal with this property. Since $Q(\mathbf{s})$ contains $j$ (which will be removed by $s_{\delta_j}$) there must exist a root of the form $\alpha_{j,\bullet}$ such that $s'_{\alpha_{j,\bullet}}\neq 0$ (in order to remove entry $j$ from $Q(\mathbf{s}')$). Since $\mathbf{s}'\in S(\varpi_i)$ this is only possible if 
$$s'_{\delta_1}=\cdots=s'_{\delta_j}=0.$$
By the minimality of $j$ we must have $s_{\delta_1}=\cdots=s_{\delta_{j-1}}=0$ and by \eqref{1234} we get the existence of $r\in\{j+1,\dots,n\}$ such that $s'_{\delta_r}=1$ but $s_{\delta_r}=0$. Again we can repeat the above argument. Since $Q(\mathbf{s}')$ contains $r$ (which will be removed by $s'_{\delta_r}$) there must exist a root of the form $\alpha_{r,\bullet}$ such that $s_{\alpha_{r,\bullet}}\neq 0$ (in order to remove entry $r$ from $Q(\mathbf{s}))$. This is a contradiction to $s_{\delta_j}\neq 0$ since $\mathbf{s}\in S(\varpi_i)$. Hence $s_{\delta_j}=s'_{\delta_j}$ for all $1\leq j\leq n$.
This implies $Q(\mathbf{s})=Q(\mathbf{s}')$ and a straightforward calculation shows now $\mathbf{s}=\mathbf{s}'$, which is a contradiction.
\end{proof}
\end{lem}
The proof of $|S(\varpi_i)|=\dim V(\varpi_i)$ is now finished as follows. By Lemma~\ref{inj}, we have $|S(\varpi_i)|\leq |KT(\epsilon\varpi_i)|$ and by \cite[Proposition 4.2]{KT90} we have that $|KT(\epsilon\varpi_i)|$ is equal to the dimension of the irreducible $\mathfrak{so}_{2n+1}$ representation of highest weight $\epsilon\varpi_i$. Now, using \cite[Theorem 2.1]{RS82} we get
$\dim V(\varpi_i)=\dim V_{\mathfrak{so}_{2n+1}}(\epsilon\varpi_i)$. Putting all together, we get
$$|S(\varpi_i)|\leq  |KT(\epsilon\varpi_i)|=\dim V_{\mathfrak{so}_{2n+1}}(\epsilon\varpi_i)=\dim V(\varpi_i).$$
The converse direction follows from the spanning property.

\section{Appendix: The exceptional cases}\label{section-appendix}
Here we will consider the exceptional cases $D(2,1;\alpha)$, $F(4)$ and $G(3)$ and construct a PBW basis parametrized by the lattice points of a polytope for a class of dominant integral typical weights. We will need a combinatorial result first whose proof is straightforward and will be omitted.
\begin{lem}\label{comb}The cardinality of the set
$$\{(a_0,a_1,\dots,a_\ell)\in \mathbb{Z}_+\times \{0,1\}^{\ell}: \sum_{i=0}^\ell a_i\leq m\}$$ is given by $2^{\ell-1}(2m-\ell+2)$ provided that $m\geq \ell-1$.
\hfill\qed
\end{lem}
We define $$P_{d}^+=\{\lambda\in P^+: \lambda \text{ typical, } k\geq | R^+_{\bar{1}}|-1\},$$
where $k$ is the number in \cite[Table 2]{K78}. The fact that $\lambda\in P^+$ is typical gives already some restrictions on $k$, e.g. in $D(2,1;\alpha)$ we must have $k\geq 2$ and in $F(4)$ we must have $k\geq 4$ and in $G(3)$ we get $k\geq 3$. Note that $k\geq | R_{\bar{1}}|-1$ is slightly stronger, e.g. in $D(2,1;\alpha)$ we have four positive odd roots. 
\subsection{}Consider the Lie superalgebra $D(2,1;\alpha)$ with distinguished simple system 
$$\{\alpha_1=\epsilon_1-\epsilon_2-\epsilon_3, \alpha_2=2\epsilon_2,\alpha_3=\epsilon_3\},$$
where $\{\epsilon_1,\epsilon_2,\epsilon_3\}$ is an orthogonal basis such that $(\epsilon_1,\epsilon_1)=\frac{-(1+\alpha)}{2}$, $(\epsilon_2,\epsilon_2)=\frac{1}{2}$, $(\epsilon_3,\epsilon_3)=\frac{\alpha}{2}.$
The set of positive roots is given by
$$R_{\bar{0}}^+=\{2\epsilon_1,2\epsilon_2,2\epsilon_3\},\ \ R_{\bar{1}}^+=\{\epsilon_1\pm \epsilon_2\pm \epsilon_3\}.$$
We fix a dominant integral typical weight $\lambda\in P_d^+$. Since $\lie g_{\bar{0}}\cong \mathfrak{sl}_2\oplus \mathfrak{sl}_2\oplus\mathfrak{sl}_2$ we must have
$$\lambda=m_1\varpi^1_1+m_2\varpi^2_1+m_3\varpi^3_1,\ \ m_1,m_2,m_2\in\mathbb{Z}_+$$ with some further restrictions on $m_1,m_2,m_3$; see for example \cite[pg. 622]{K78} (in this case we have $m_1=k$). We have with Lemma~\ref{dimzer}
$$\text{dim } V(\lambda)=16 (m_1-1)(m_2+1)(m_3+1).$$
By the PBW theorem it is clear that the following elements span $\gr V(\lambda)$:
\begin{equation}\label{span}\left(\prod_{\beta\in R_{\bar{1}}^+}x_{-\beta}^{s_{\beta}}\right)x^{s_{2\epsilon_1}}_{-2\epsilon_1}x^{s_{2\epsilon_2}}_{-2\epsilon_2}x^{s_{2\epsilon_3}}_{-2\epsilon_3},\ \  s_{2\epsilon_1}\leq m_1,\ s_{2\epsilon_2}\leq m_2,\ s_{2\epsilon_3}\leq m_3,\ s_{\beta}\in\{0,1\}.\end{equation}
We will need one further relation in \eqref{span} in order to obtain a basis. Recall that $\gamma=2\alpha_1+\alpha_2+\alpha_3$ and hence $x_{-2\epsilon_1}^{m_1+1}v_{\lambda}=0$ in $\gr V(\lambda)$. By applying the operators $\partial_{\beta},\beta\in R_{\bar{1}}^+$ we obtain:
\begin{align*}0&=\partial_{\epsilon_1-\epsilon_2-\epsilon_3}\partial_{\epsilon_1-\epsilon_2+\epsilon_3}\partial_{\epsilon_1+\epsilon_2-\epsilon_3}\partial_{\epsilon_1+\epsilon_2+\epsilon_3}x_{-2\epsilon_1}^{m_1+1}&\\&=x_{-2\epsilon_1}^{m_1-3}x_{-\epsilon_1-\epsilon_2-\epsilon_3}x_{-\epsilon_1-\epsilon_2+\epsilon_3}x_{-\epsilon_1+\epsilon_2-\epsilon_3}x_{-\epsilon_1+\epsilon_2+\epsilon_3}+x_{-2\epsilon_1}^{m_1-2}x_{-2\epsilon_2}x_{-\epsilon_1+\epsilon_2-\epsilon_3}x_{-\epsilon_1+\epsilon_2+\epsilon_3}&\\&\hspace{0,5cm}+x_{-2\epsilon_1}^{m_1-2}x_{-\epsilon_1-\epsilon_2+\epsilon_3}x_{-2\epsilon_3}x_{-\epsilon_1+\epsilon_2+\epsilon_3}.\end{align*}
Now choosing an appropriate order on the set of positive roots, we can impose the following relation in \eqref{span}
$$\sum_{\beta\in R_{\bar{1}}^+} s_{\beta}+s_{2\epsilon_1}\leq m_1.$$
We claim that the set $\{x^{\mathbf{s}}v_{\lambda}: \mathbf{s}\in S(\lambda)\}$
is a basis of $\gr V(\lambda)$, where 
$$S(\lambda)=\Big\{\mathbf{s}: \sum_{\beta\in R_{\bar{1}}^+} s_{\beta}+s_{2\epsilon_1}\leq m_1,\ s_{2\epsilon_2}\leq m_2,\ s_{2\epsilon_3}\leq m_3,\ \forall \beta\in R_{\bar{1}}^+:  \ s_{\beta}\leq 1\Big\}.$$
This follows immediately if we can show that the cardinality of $S(\lambda)$ equals the dimension of $V(\lambda)$. But this is clear with Lemma~\ref{comb}.
\subsection{}Here we consider the Lie superalgebra $F(4)$ with distinguished simple system 
$$\{\alpha_1=\frac{1}{2}(\delta-\epsilon_1-\epsilon_2-\epsilon_3), \alpha_2=\epsilon_3,\alpha_3=\epsilon_2-\epsilon_3,\alpha_4=\epsilon_1-\epsilon_2\}$$
and positive roots
$$R_{\bar{0}}^+=\{\delta,\epsilon_1\pm\epsilon_2, \epsilon_2\pm\epsilon_3,\epsilon_1\pm\epsilon_3, \epsilon_1,\epsilon_2,\epsilon_3\},\ \ R_{\bar{1}}^+=\{\frac{1}{2}(\delta\pm\epsilon_1\pm\epsilon_2\pm\epsilon_3)\},$$
where  $(\epsilon_i,\epsilon_j)=\delta_{i,j},\ (\epsilon_i,\delta)=0, \ (\delta,\delta)=-3$. Again we fix a dominant integral weight $\lambda\in P_d^+$. Since $\lie g_{\bar{0}}\cong \mathfrak{sl}_2\oplus \mathfrak{so}_7$ we must have 
$$\lambda=m_1\varpi^1_1+k_1\varpi^2_1+k_2\varpi^2_2+k_3\varpi^2_3,\ \ m_1,k_1,k_2,k_3\in\mathbb{Z}_+$$
with some further restrictions on $m_1,k_1,k_2,k_3$; see for example \cite[pg. 622]{K78} (in this case we have $m_1=k$). Setting $\lambda'=\lambda-m_1\varpi_1^1$ we get with Lemma~\ref{dimzer} that
$$\text{dim } V(\lambda)=2^8 (m_1-3)\text{ dim } V_{\mathfrak{so}_7}(\lambda').$$
Recall that a convex polytope parametrizing a basis of $\gr V_{\mathfrak{so}_7}(\lambda')$ has been determined in \cite[Theorem 5.2]{BK19}; we denote this polytope by $S_{\mathfrak{so}_7}(\lambda')$. In fact, we can choose any basis but we favor a basis parametrized by the lattice points of a polytope.
A similar calculation as in the previous subsection and Proposition~\ref{reduction} shows that the following set is a basis of $\gr V(\lambda)$:
\begin{equation}\label{123}\{x^{\mathbf{s}}v_{\lambda}: \mathbf{s}\in S(\lambda)\},\end{equation}
where 
$$S(\lambda)=\Big\{\mathbf{s}: \sum_{\beta\in R_{\bar{1}}^+} s_{\beta}+s_{\delta}\leq m_1, \ \ (s_{\alpha})_{\alpha\in R_{\bar{0}}^+\backslash\{\delta\}}\in S_{\mathfrak{so}_7}(\lambda')\Big\}.$$
The fact that \eqref{123} is a spanning set is done in a similar fashion by applying suitable differential operators and choosing an appropriate order on the set of positive roots. The fact that the cardinality of $S(\lambda)$ is equal to $\text{dim } V(\lambda)$ follows again by Lemma~\ref{comb}.
\subsection{}Here we consider the Lie superalgebra $G(3)$ with distinguished simple system 
$$\{\alpha_1=\delta+\epsilon_3, \alpha_2=\epsilon_1,\alpha_3=\epsilon_2-\epsilon_1\}.$$
The set of positive roots is given by
$$R_{\bar{0}}^+=\{2\delta,\epsilon_1,\epsilon_2-\epsilon_1,\epsilon_2,\epsilon_1+\epsilon_2,2\epsilon_1+\epsilon_2,\epsilon_1+2\epsilon_2\},\ \ R_{\bar{1}}^+=\{\delta,\ \delta\pm\epsilon_i: 1\leq i\leq 3\},$$
where $\epsilon_1+\epsilon_2+\epsilon_3=0$ and $(\epsilon_i,\epsilon_j)=1-3\delta_{i,j},\ (\delta,\delta)=2,\ (\epsilon_i,\delta)=0.$ Again we fix a dominant integral weight $\lambda\in P_d^+$. Since $\lie g_{\bar{0}}\cong \mathfrak{sl}_2\oplus G_2$ we must have 
$$\lambda=m_1\varpi^1_1+k_1\varpi^2_1+k_2\varpi^2_2,\ \ m_1,k_1,k_2\in\mathbb{Z}_+$$
with some further restrictions on $m_1,k_1,k_2$; see for example \cite[pg. 633]{K78} (in this case we have $m_1=k$). Setting $\lambda'=\lambda-m_1\varpi_1^1$ we get with Lemma~\ref{dimzer} that
$$\text{dim } V(\lambda)=64 (2m_1-5) \text{ dim }V_{G_2}(\lambda').$$
Recall that a convex polytope parametrizing a basis of $\gr V_{G_2}(\lambda')$ has been determined in \cite[Theorem 1]{G15}; we denote this polytope by $S_{G_2}(\lambda')$. Again for our purposes we could choose any basis. A similar calculation as in the previous subsections and Proposition~\ref{reduction} shows that the following set is a basis of $\gr V(\lambda)$:
\begin{equation}\label{12345}\{x^{\mathbf{s}}v_{\lambda}: \mathbf{s}\in S(\lambda)\},\end{equation}
where 
$$S(\lambda)=\Big\{\mathbf{s}: \sum_{\beta\in R_{\bar{1}}^+} s_{\beta}+s_{2\delta}\leq m_1,\ s_{\beta}\leq 1 \ \forall \beta\in R_{\bar{1}}^+, \ (s_\tau)_{\tau\in R_{\bar{0}}^+\backslash\{2\delta\}}\in S_{G_2}(\lambda')\Big\}.$$
The fact that \eqref{12345} is a spanning set is done by applying differential operators corresponding to positive odd roots to the element $x_{-2\delta}^{m_1+1}$ (recall that this element vanishes in  $\gr V(\lambda)$). The fact that the cardinality of $S(\lambda)$ is equal to $\text{dim } V(\lambda)$ is done similarly by using Lemma~\ref{comb}.
%

\bibliographystyle{plain}
\bibliography{bibfile}
\end{document}